\documentclass{article}
\usepackage{epsfig}

\newtheorem{theorem}{Theorem}[section]
\newtheorem{corollary}[theorem]{Corollary}

\begin{document}

\title{Generalized persistency of excitation}

\author{Sergey Nikitin \\
Department of Mathematics\\
and\\
Statistics\\
Arizona State University\\
Tempe, AZ 85287-1804}

\maketitle

{\bf Abstract} The paper presents the generalized persistency of excitation conditions. They are not only valid for much broader range of applications than their classical counterparts but also elegantly prove the validity of the latter. The novelty and the significance of the approach presented in this publication is due to employing the time averaging technique.

\vspace{0.1cm}

\section{Introduction}

The persistency of excitation conditions appear in numerous applications related to system identification, learning, adaptation, parameter estimation. They guarantee the convergence of the adaptation procedures based on the ideas of gradient and least-squares algorithms. An introduction into this topic can be found, e.g., in chapter 2 of the book \cite{Sastry}. The classical version of the persistency of excitation conditions can be characterized in terms of the asymptotic stability  of the linear system
\begin{equation}
\label{system}
\dot x = - P(t) x,\;\;\;\;P(t)=P(t)^T \ge 0,\;\;\forall\;t\ge0.
\end{equation}
Namely, the linear system is uniformly asymptotically stable if $P(t)$ is persistently exciting, i.e., there exist positive real numbers $\alpha,\;\;\beta,\;\;\delta$ such that
\begin{equation}
\label{classical1}
\beta \cdot I \ge \int_{t}^{t+\delta } P(\tau) d\tau \ge \alpha \cdot I \;\; \forall t \ge 0,
\end{equation}
where $I$ is the identity matrix. The detailed analysis of the conditions (\ref{classical1}) can be found in many publications (see, e.g., \cite{Anderson}, \cite{Janecki}, \cite{Loria},  \cite{Morgan1}, \cite{Morgan2}, \cite{Kreisselmeier}, \cite{Sondhi},  \cite{Sastry}). 

We address (\ref{classical1}) as classical persistency of excitation conditions. They impose the restrictions that are uniform in time. Moreover, they tacitly demand the exponential convergence of the corresponding adaptation procedures.  On the other hand, due to wide range of applications the classical conditions might be a burden for solutions of important problems. This publication presents the new generalized persistency of excitation conditions that do not impose any unnatural (uniform in time, exponential convergence) restrictions. Moreover, the classical version easily follows from the new generalized conditions. 

\vspace{0.1cm}

 In order to create the generalized version for the persistency of excitation this paper uses the approach similar in the spirit to the time-averaging developed in \cite{Aeyels} and the technique originated in the theory of linear boundary value problems for partial differential equations (see , e.g., \cite{Triebel}). In other words, we introduce a generalized definition of a solution for an ordinary differential equation. That definition follows the widely accepted ideology of distributions. Then we formulate our new necessary and sufficient conditions for a time-varying system (\ref{system}) to be asymptotically stable. Finally, we formulate corollaries of the main result and present an example illustrating the generalized persistency of excitation conditions.

\vspace{0.1cm}

\section{Preliminaries}
Consider a system
\begin{equation}
\label{linearSystem}
\dot x =  -P(t)x
\end{equation}
where $x\in {\rm R}^n,$ $ {\rm R}^n$ -- $n$-dimensional linear real space. $P(t)$ is a time-dependent  matrix such that
$$
P(t)=P(t)^{T} \ge 0 \;\; \forall \; t \in {\rm R},
$$
where the inequality $P(t) \ge 0$ is understood in the following sense.
Given two $n\times n$ symmetric matrices $A$ and $B$ we write
$$
A \ge B
$$
if
$$
\langle x, A x \rangle  \ge \langle x, B x \rangle \;\;\;\forall \;x \;\in {\rm R}^n.
$$
Throughout the paper we assume that  ${\rm R}^n$ is equipped with the scalar product and $\| x \|$ denotes the magnitude of $x,$ i.e
$\| x \| =\sqrt {\langle x , x \rangle },$ where $\langle x , x \rangle$ is the scalar product of $x$ with itself. 

\vspace{0.1cm}

We assume that $P(t)$ is a time-dependent $L_{1, loc}$-matrix in the following sense. For any real numbers $b>a$ and for any $x,\;y\in  {\rm R}^n $ we have
$$
\int_a^b \mid \langle y, P(t) x \rangle \mid dt < \infty.
$$

Consider the initial value problem
\begin{eqnarray}
\label{initialValueProblem}
\dot x(t) &=& -  P(t)x(t) ,       \nonumber\\
&&\\
x(0)&=&x_0, \nonumber
\end{eqnarray}
where $P(t) \in L_{1, loc}.$
Its solution is defined to be an $L_{1, loc}$ vector function $x(t)$
such that for any infinitely differentiable function $\varphi (t) \in C^{\infty }$ (both $x(t)$ and $\varphi (t)$ take its values from ${\rm R}^n$) with compact support (that means $\varphi (t)=0$ outside an interval from ${\rm R})$ we have
\begin{equation}
\label{generilizedSolution}
- \int_{0}^\infty \langle\frac{d}{dt} \varphi (t), x(t) \rangle dt = \langle \varphi (0),  x_0 \rangle - \int_{0}^\infty \langle \varphi (t), P(t)x(t) \rangle dt .
\end{equation}

It is well-known (see e.g., \cite{Coddington},  \cite{Sansone}) that the solution $x(t,x_0)$ for (\ref{initialValueProblem}) exists and unique. Moreover, it is a continuous function of time. Indeed, consider the Picard's sequence
$$
y_n(t) = x_0 - \int_{0}^t P(\tau) y_{n-1}(\tau) d\tau, 
$$
where $y_0(t)=x_0$ and $n=1,\;2,\;\dots .$ Since $P(t) \in L_{1, loc}$ the Picard's sequence $\{y_n(t)\}_n$ converges (point-wise) to a continuous function $y(t).$ Let us show that $y(t)$ satisfies (\ref{generilizedSolution}). Integrating by parts
$$
  -\int_{0}^\infty \langle\frac{d}{dt} \varphi (t), y_n(t) \rangle dt
$$
we obtain
$$
 -\int_{0}^\infty \langle\frac{d}{dt} \varphi (t), y_n(t) \rangle dt =  \langle \varphi (0), x_0 \rangle - \int_{0}^\infty  \langle \varphi (t),  P(t) y_{n-1}(t) \rangle d t
$$
Since both $ y_n(t)$ and $y_{n-1}(t)$ converge to $y(t)$ as $n\to \infty$ we arrive at
$$
-\int_{0}^\infty \langle\frac{d}{dt} \varphi (t), y(t) \rangle dt =  \langle \varphi (t_0), x_0 \rangle - \int_{0}^\infty  \langle \varphi (t),  P(t), y(t) \rangle d t.
$$
Hence, $y(t)=x(t,x_0)$ is the solution for (\ref{initialValueProblem}) in the sense (\ref{generilizedSolution}). 

The goal of this paper is to find necessary and sufficient conditions for  the solution $x(t,x_0)$ of the  system (\ref{initialValueProblem}) to satisfy:
$$
\lim_{t \to \infty } x(t,x_0) = 0  \;\;\forall\; x_0 \in {\rm R}^n ,  
$$
and the equilibrium $x=0$ is stable.

\section{Necessary and sufficient conditions}

Our main goal is to study asymptotic stability of the origin for the system  (\ref{linearSystem}). Consider a  real positive number $S$ and a continuous real function $\omega_S (t) $ such that 
$$
\omega_S (t) > 0 \;\;\mbox{ for } \;\;0\le t < S 
$$
and
$$
\omega_S (t) = 0 \;\;\mbox{ for } \;\;t \ge S.
$$
We also assume that $\omega_S (t)$ is differentiable almost everywhere on ${\rm R}.$ For the sake of brevity, we address $\omega_S (t)$ as a truncation function in the sequel.
If the origin for the system (\ref{linearSystem}) is asymptotically stable then 
$$
\lim_{t \to \infty } \| x(t,x_0) \|^2 =0 \;\;\forall\; x_0 \in {\rm R}^n .
$$
Hence, for any fixed real positive number $S$ we have
$$
\lim_{t \to \infty } \int_0^S \omega_S (\tau)\| x(t+\tau,x_0) \|^2 d\tau =0 \;\;\forall\; x_0 \in {\rm R}^n .
$$
Consider more closely the integral
$$
\int_0^S \omega_S (\tau) \| x(t+\tau,x_0) \|^2 d\tau.
$$
It follows from
$$
\frac{d}{dt} \| x(t,x_0) \|^2 = - 2\cdot \langle x(t,x_0), P(t) x(t,x_0) \rangle \le 0
$$
that $\| x(t_1,x_0) \|^2 \ge \| x(t_2,x_0) \|^2$  for $t_2 > t_1 .$
That means
$$
\int_0^S \omega_S (\tau) \cdot \| x(t+\tau,x_0) \|^2 d\tau \ge \| x(t+S ,x_0) \|^2 \cdot \int_0^S \omega_S (\tau) d\tau .
$$
If we find conditions that guarantee 
$$
\int_0^S  \omega_S (\tau) \| x(t+\tau,x_0) \|^2 d\tau \;\;\to \;\; 0 \;\;\mbox{ as }\;\;t\;\;\to \;\; \infty 
$$
then that will imply the asymptotic stability of the origin for the system (\ref{linearSystem}).

Differentiating the integral 
$$
\int_0^S  \omega_S (\tau) \cdot \| x(t+\tau,x_0) \|^2 d\tau
$$
with respect to time yields
$$
\frac{d}{dt} \int_0^S \omega_S (\tau) \cdot \| x(t+\tau,x_0) \|^2 d\tau = - 2 \cdot \int_0^S \omega_S (\tau) \cdot \langle x(t+\tau,x_0), P(t + \tau)  x(t+\tau,x_0) \rangle d\tau 
$$
Replacing $P(t + \tau)$ with
$$
\frac{d}{d\tau }\int_0^\tau P(t + \theta )d\theta
$$
and integrating by parts leads us to the following important formula
\begin{equation}
\label{main_form}
\frac{d}{dt} \int_0^S \omega_S (\tau) \cdot \| x(t+\tau,x_0) \|^2 d\tau = - 2\int_0^S  \langle x(t+\tau,x_0), A(t,\tau)  x(t+\tau,x_0) \rangle d\tau,
\end{equation}
where
$$
A(t,\tau) = - \frac{d}{d\tau } \omega_S (\tau) \cdot \int_0^\tau P(t+ \theta ) d\theta + \omega_S (\tau) \cdot \frac{d}{d\tau}(\int_0^\tau P(t+ \theta ) d\theta )^2.  
$$
Let 
$$
\lambda_{min} (t,\tau),  \;\lambda_{max} (t,\tau)
$$
denote minimal and maximal eigenvalues of $A(t,\tau).$ Then the next theorem gives us necessary and sufficient conditions for the  system (\ref{linearSystem}) to be asymptotically stable at the origin. Notice that one is assured by $P(t) \in L_{1,loc}$ that the integral expressions in the next theorem are well defined.
\begin{theorem} ({\bf Generalized persistency of excitation})
\label{GPE}
If there exist a real number $S > 0 $ and a truncation function $\omega_S (t)$  such that
$$
\lim_{t \to \infty } \sup \int_0^t \int_0^S \lambda_{min} (\nu,\tau) d\tau d\nu  = \infty
$$ 
then the system (\ref{linearSystem}) is asymptotically stable at the origin.
On the other hand, if 
$$
\;\;\;\lim_{t \to \infty } \inf \int_0^t \int_0^S \lambda_{max} (\nu,\tau) d\tau d\nu  < \infty 
$$ 
then the system (\ref{linearSystem}) is not asymptotically stable at the origin. 
\end{theorem}
{\bf Proof.}
It follows from (\ref{main_form}) that
\begin{eqnarray*}
\int_0^S \omega_S(\tau ) \cdot \| x(t+\tau,x_0) \|^2 d\tau &\le& -\int_0^t \int_0^S \lambda_{min} (\nu,\tau) \| x(\nu+\tau,x_0) \|^2 d\tau d\nu \\
&&  + \int_0^S \omega_S(\tau ) \| x(\tau,x_0) \|^2 d\tau
\end{eqnarray*}
Due to monotonicity of $\| x(t,x_0) \|^2$ we have
\begin{eqnarray*}
(\int_0^S \omega_S(\tau ) d\tau ) \cdot  \| x(t+S,x_0) \|^2 &\le& -\int_0^t (\int_0^S \lambda_{min} (\nu,\tau)d\tau ) \| x(\nu+S,x_0) \|^2 d\nu\\
 &&+(\int_0^S \omega_S(\tau ) d\tau ) \cdot  \| x_0 \|^2.
\end{eqnarray*}
It follows from the very well known Gronwall inequality ( see, e.g., \cite{Rouche}) that
$$
\| x(t+S,x_0) \|^2 \le \| x_0 \|^2 \cdot \exp{\{- \frac{\int_0^t \int_0^S \lambda_{min} (\nu,\tau)d\tau d\nu } {\int_0^S \omega_S(\tau ) d\tau } \}}
$$
Thus,
$$
\lim_{t \to \infty } \sup \int_0^t \int_0^S \lambda_{min} (\nu,\tau) d\tau d\nu  = \infty
$$ 
implies that the system (\ref{linearSystem}) is asymptotically stable.

On the other hand,  (\ref{main_form}) leads us to
\begin{eqnarray*}
\int_0^S\omega_S(\tau ) \cdot \| x(t+\tau,x_0) \|^2 d\tau &\ge&  -\int_0^t \int_0^S \lambda_{max} (\nu,\tau) \| x(\nu+\tau,x_0) \|^2 d\tau d\nu \\
&& + \int_0^S \omega_S(\tau ) \cdot \| x(\tau,x_0) \|^2 d\tau
\end{eqnarray*}
Since  the derivative  $\frac{d}{dt}\| x(t,x_0) \|^2$   is not positive we have $\| x(\nu+\tau,x_0) \|^2 \le \| x(\nu,x_0) \|^2 $ and 
\begin{eqnarray*}
(\int_0^S\omega_S(\tau )d\tau) \cdot \| x(t,x_0) \|^2 &\ge & -\int_0^t (\int_0^S \lambda_{max} (\nu,\tau)d\tau) \|  x(\nu,x_0)\|^2 d\nu \\
&& +(\int_0^S\omega_S(\tau )d\tau)  \| x(S,x_0)  \|^2
\end{eqnarray*}
for $t \ge S .$ After solving this inequality we obtain
$$
\| x(t,x_0) \|^2 \ge \| x(S,x_0) \|^2 \cdot \exp{\{- \frac{\int_0^t \int_0^S \lambda_{max} (\nu,\tau)d\tau d\nu}{\int_0^S\omega_S(\tau )d\tau} \}}
$$
for $t \ge S .$
Hence, if 
$$
\;\;\;\lim_{t \to \infty } \inf \int_0^t \int_0^\S \lambda_{max} (\nu,\tau) d\tau d\nu  < \infty 
$$ 
then the system (\ref{linearSystem}) is not asymptotically stable.
{\bf Q.E.D.}

Notice that we owe the success in proving Theorem  \ref{GPE} to the new idea that suggests to consider 
$$
 \int_0^S \omega_S(\tau ) \cdot \| x(t+\tau,x_0) \|^2 d \tau
$$
instead of $ \| x(t,x_0) \|^2.$ This approach seems to have further important consequences not only for control theory but also for studies of general dynamical systems. 

\vspace{0.1cm}

After integrating by parts we have
\begin{equation}
\label{effective}
\int_0^S A(t,\tau) d\tau =  -\int_0^S \frac{d}{d\tau } \omega_S (\tau) \cdot \{\int_0^\tau P(t+ \theta ) d\theta + (\int_0^\tau P(t+ \theta ) d\theta)^2 \}d\tau.
\end{equation}
Let
$$
\gamma_{min} (t,\tau),  \;\gamma_{max} (t,\tau)
$$
denote minimal and maximal eigenvalues of 
$$
\int_0^\tau P(t+ \theta ) d\theta.
$$
 Then we can reformulate Theorem \ref{GPE} as follows. 
\begin{corollary}
If there exist a real number $S > 0 $ and a truncation function $\omega_S (t)$  such that 
$$
\frac{d}{d\tau } \omega_S (\tau) \le 0 \;\;\; \mbox{ for } \tau < S
$$
and
$$
\lim_{t \to \infty } \sup \{-\int_0^t \int_0^S \frac{d}{d\tau } \omega_S (\tau) (\gamma_{min} (\nu,\tau) + (\gamma_{min} (\nu,\tau))^2) d\tau d\nu\}  = \infty
$$
then the system (\ref{linearSystem}) is asymptotically stable at the origin.
On the other hand, if
$$
\lim_{t \to \infty } \inf \{-\int_0^t \int_0^S \frac{d}{d\tau } \omega_S (\tau) (\gamma_{max} (\nu,\tau) + (\gamma_{max} (\nu,\tau))^2) d\tau d\nu\}  < \infty
$$
then the system (\ref{linearSystem}) is not asymptotically stable at the origin.
\end{corollary}
{\bf Proof.}  
Consider the unit eigenvectors $\psi_{min}(\nu,\tau),\;\;\psi_{max}(\nu,\tau) $ corresponding to minimal and maximal eigenvalues of $\int_0^\tau P(\nu + \theta ) d\theta,$ 
$$
\gamma_{min} (\nu,\tau) = \langle \psi_{min}, (\int_0^\tau P(\nu+ \theta ) d\theta) \psi_{min}\rangle,\;\;\; \langle\psi_{min} (\nu,\tau) , \psi_{min} (\nu,\tau) \rangle = 1,
$$
and
$$
\gamma_{max} (\nu,\tau) = \langle \psi_{max}, (\int_0^\tau P(\nu+ \theta ) d\theta) \psi_{max}\rangle,\;\;\; \langle\psi_{max} (\nu,\tau) , \psi_{max} (\nu,\tau) \rangle = 1.
$$
Then, taking into account that
$$
\langle \frac{d}{d\tau} \psi_{min},\psi_{min}\rangle =0,\;\; \langle \frac{d}{d\tau} \psi_{max},\psi_{max}\rangle =0
$$ 
we obtain
$$
 \langle \psi_{min},\frac{d}{d\tau }  (\int_0^\tau P(\nu+ \theta ) d\theta)^2  \psi_{min} \rangle =\frac{d}{d\tau } \langle \psi_{min},  (\int_0^\tau P(\nu+ \theta ) d\theta)^2  \psi_{min} \rangle
$$
and
$$
 \langle \psi_{max},\frac{d}{d\tau }  (\int_0^\tau P(\nu+ \theta ) d\theta)^2  \psi_{max} \rangle =\frac{d}{d\tau } \langle \psi_{max},  (\int_0^\tau P(\nu+ \theta ) d\theta)^2  \psi_{max} \rangle.
$$
Consequently,
\begin{equation}
\label{min}
\langle \psi_{min} , A(\nu,\tau)  \psi_{min} \rangle = -\frac{d}{d\tau } \omega_S (\tau) \gamma_{min} (\nu,\tau) + \omega_S (\tau) \frac{d}{d\tau }(\gamma_{min} (\nu,\tau))^2 
\end{equation}
and
\begin{equation}
\label{max}
\langle \psi_{max} , A(\nu,\tau)  \psi_{max} \rangle = -\frac{d}{d\tau } \omega_S (\tau) \gamma_{max} (\nu,\tau) +\omega_S (\tau) \frac{d}{d\tau }(\gamma_{max} (\nu,\tau))^2 .
\end{equation}
Let $\lambda_{min}(\nu,\tau)$ and $\lambda_{max}(\nu,\tau)$ denote minimal and maximal eigenvalues for $A(\nu,\tau) ,$ respectively. Then, after integrating by parts
$$
\int_0^S \lambda_{min}(\nu,\tau)d\tau = -\int_0^S \frac{d}{d\tau } \omega_S (\tau) (\gamma_{min} (\nu,\tau) + (\gamma_{min} (\nu,\tau))^2) d\tau
$$
follows from (\ref{min}) and
$$
\int_0^S \lambda_{max}(\nu,\tau)d\tau = - \int_0^S \frac{d}{d\tau } \omega_S (\tau) (\gamma_{max} (\nu,\tau) + (\gamma_{max} (\nu,\tau))^2) d\tau
$$
follows from (\ref{max}).
{\bf Q.E.D.}

\vspace{0.1cm}

If
$$
\omega_S(\tau ) =  \left\{
                      \begin{array}{cc}
                          (S - \tau ) \; & \mbox{ for } \tau < S  \\
                          0 \;\;\;\;\;\;\;\;\;\;\;\;\; & \mbox{ for } \tau \ge S
                      \end{array}
                   \right. 
$$
then we obtain the following important corollary of Theorem  \ref{GPE}.
\begin{corollary} 
\label{generalAppl} 
If there exists a real number $S > 0 $ such that
$$
\lim_{t \to \infty } \sup \int_0^t \int_0^S \gamma_{min} (\nu,\tau) d\tau d\nu  = \infty
$$
then the system (\ref{linearSystem}) is asymptotically stable at the origin.
\end{corollary}

Though Corollary \ref{generalAppl} gives us only a sufficient condition of asymptotic stability for the system (\ref{linearSystem}) its simple form makes it valuable for practical applications. 

\vspace{0.1cm}

At the conclusion of this section we present a simple and  elegant proof for the classical persistency of excitation conditions. 
\begin{corollary} (Classical Persistency of Excitation)
If there exist real numbers $\alpha >0,\;\;\delta > 0$   such that 
\begin{equation}
\label{classic}
\int_{0}^\delta  P(t+s)  ds \ge \alpha I\;\;\;\forall\;t \ge 0
\end{equation}
then the system (\ref{linearSystem}) is asymptotically stable.
\end{corollary}
{\bf Proof.} It follows from (\ref{classic}) that for the minimal eigenvalue $\gamma_{min}(\nu,\tau)$ from Corollary \ref{generalAppl} we have 
$$
\gamma_{min}(\nu,\tau) \ge \alpha \;\;\;\mbox{ for } \tau \ge \delta.
$$
Hence, if we take $S>\delta$ then 
$$
\int_0^t \int_0^S \gamma_{min} (\nu,\tau) d\tau d\nu \ge \int_0^t \int_\delta^S \gamma_{min} (\nu,\tau) d\tau d\nu \ge (S-\delta) t \alpha 
$$
and
$$
\lim_{t \to \infty } \sup \int_0^t \int_0^S \gamma_{min} (\nu,\tau) d\tau d\nu \ge \lim_{t \to \infty } (S-\delta) t  \alpha = \infty.
$$
\vspace{0.1cm}

{\bf Q.E.D.}

\section{Example}

Theorem \ref{GPE} and its corollaries find many important applications. This section illustrates how one can use them in order to verify the persistency of excitation conditions.

\vspace{0.1cm}

Consider the system
\begin{equation}
\label{ex1}
\dot x = -\left( \begin{array}{cc}
                    \Xi(t) & 0  \\
                      0 &  1 - \Xi(t)
                \end{array}
         \right) x  
\end{equation}
where $\Xi(t)$ is a characteristic function of a closed subset $C \subset {\rm R},$
$$
\Xi(t) = \left\{ \begin{array}{ccc}
                1 & \mbox{ for } & t \in C,\\
                0 &  \mbox{otherwise}& 
                 \end{array}  
         \right.          
$$
Consider the sequence of real numbers $\{a_n\}_{n=0}^\infty$ defined as
$$
a_{n+1} = (n + 1) + a_n,\;\;\;
$$
where $a_0 = 0 .$ If we define the closed subset $C$ as 
$$
C= \{ t \in {\rm R};\;\;\exists \;n > 1\;\;\mbox{ such that } a_n \le t \le a_n+1 \}
$$
then the classical condition (\ref{classical1}) is not valid. However, the system is persistently exciting due to Corollary \ref{generalAppl}. Indeed, under the assumption that $0< \tau <1$ we have
$$
\gamma_{min}(\nu, \tau )=\left\{ \begin{array}{ccc}
                \tau - a_n+ \nu & \exists \; n >1 & a_n - \tau \le \nu \le a_n - \frac{\tau}2 , \\
                a_n - \nu  &  \exists \; n >1 & a_n - \frac{\tau}2  \le \nu \le a_n ,\\
                \tau - (a_n+1)+ \nu &  \exists \; n >1 & a_n+1 - \tau \le \nu \le a_n +1- \frac{\tau}2 , \\
                a_n+1 - \nu  &  \exists \; n >1 & a_n+1 - \frac{\tau}2  \le \nu \le a_n+1,\\
                    0 &  \mbox{otherwise}&
                 \end{array}  
         \right.
$$
Making use of Fubini theorem we obtain
$$
\int_0^t \int_0^S \gamma_{min} (\nu,\tau) d\tau d\nu = \int_0^S \int_0^t \gamma_{min} (\nu,\tau) d\nu d\tau .
$$
If we take $0< S <1 $ and  
$$
  a_n + 1  \le t \le a_{n+1} - 1 \mbox{ for some } n > 1,
$$ 
then
$$
\int_0^S \int_0^t \gamma_{min} (\nu,\tau) d\nu d\tau = (n-1)\frac{S^3}6.
$$
Hence, Corollary \ref{generalAppl} implies that the system is persistently exciting.

\end{document}